\documentclass[12pt]{amsart}
\font\emailfont=cmtt10

\usepackage{amsmath,amsthm,amsfonts,amscd,flafter,epsf}

\newcommand{\HF}{HF}

\newtheorem{theorem}{Theorem}[section]

\def\endproof{\relax\ifmmode\expandafter\endproofmath\else
  \unskip\nobreak\hfil\penalty50\hskip.75em\hbox{}\nobreak\hfil\bull
  {\parfillskip=0pt \finalhyphendemerits=0 \bigbreak}\fi}
\def\endproofmath$${\eqno\bull$$\bigbreak}
\def\bull{\vbox{\hrule\hbox{\vrule\kern3pt\vbox{\kern6pt}\kern3pt\vrule}\hrule}}

\newcommand{\Q}{\mathbb{Q}}
\newcommand{\R}{\mathbb{R}}

\newcommand{\C}{\mathbb{C}}

\newcommand{\Z}{\mathbb{Z}}

\newcommand{\cm}{\cdot}

\newcommand{\Nbd}[1]{{\mathrm{nd}}(#1)}
\newcommand{\nbd}[1]{\Nbd{#1}}
\newcommand{\CDisk}{\mathbb D}

\newcommand{\ModSWfour}{\mathcal{M}}
\newcommand{\ModFlow}{\ModSWfour}

\newcommand{\DBar}{\overline{\partial}}

\newcommand\abuts\Rightarrow
\newcommand\Sym{\mathrm{Sym}}

\newcommand\relspinc{\underline{\spinc}}

\newcommand\x{\mathbf x}
\newcommand\w{\mathbf w}
\newcommand\z{\mathbf z}

\newcommand\y{\mathbf y}

\newcommand\ModSphere{\ModFlow\left({\mathbb S}\longrightarrow 
\Sym^{g-1}(\Sigma_{1})\times \Sym^2(\Sigma_{2})\right)}
\newcommand\ModSpheres\ModSphere
\newcommand\CF{CF}

\newcommand\CFa{\widehat{CF}}
\newcommand\CFp{\CFb}
\newcommand\CFm{\CF^-}

\newcommand\HFp{\HFb}

\newcommand\HFm{\HF^-}
\newcommand\CFinf{CF^\infty}
\newcommand\HFinf{HF^\infty}
\newcommand\CFb{CF^+}
\newcommand\HFa{\widehat{HF}}
\newcommand\HFb{HF^+}
\newcommand\gr{\mathrm{gr}}
\newcommand\Mas{\mu}
\newcommand\UnparModSp{\widehat \ModSp}
\newcommand\UnparModFlow\UnparModSp
\newcommand\Mod\ModSp

\newcommand{\cald}{{\mathcal D}}

\newcommand\PD{\mathrm{PD}}

\newcommand{\spinc}{\mathfrak s}

\newcommand\ModMaps{\mathcal M}
\newcommand\ModSp\ModMaps

\newcommand\Ta{{\mathbb T}_{\alpha}}
\newcommand\Tb{{\mathbb T}_{\beta}}

\newcommand\alphas{\mbox{\boldmath$\alpha$}}

\newcommand\betas{\mbox{\boldmath$\beta$}}

\newcommand\Dual{\mathcal D}
\newcommand\Duality\Dual

\newcommand\Da{\widehat\partial}

\newcommand\HH{\mathbb{H}}
\newcommand\lk{\mathrm{lk}}

\newcommand\orL{\vec{L}}
\newcommand\CFLa{\widehat{\mathrm CFL}}

\newcommand\HFLa{\widehat {\mathrm{HFL}}}

\newcommand\spincrel\relspinc

\newcommand\HFK{\mathrm{HFK}}

\newcommand\HFKa{\widehat\HFK}

\keywords{Heegaard diagrams, Floer homology, Thurston norm}

\title[{Heegaard diagrams and Floer homology}]
{Heegaard diagrams and Floer homology}

\author[Peter Ozsv{\'a}th]{Peter Ozsv\'ath}
\address{Department of
Mathematics, Columbia University, 
New York, NY 10027 \newline
\indent{\emailfont{petero@math.columbia.edu}}}
\thanks{PSO was supported by NSF grant number DMS-0505811}

\author[Zolt{\'a}n Szab{\'o}]{Zolt{\'a}n Szab{\'o}} 
\address{Department of
Mathematics, Princeton University, New Jersey 08544 \newline
\indent{\emailfont{szabo@math.princeton.edu}}}
\thanks{ZSz was supported by NSF grant number DMS-0406155}

\begin{document}

\begin{abstract}  
We review the construction of Heegaard Floer homology for closed
three-manifolds and also for knots and links in the three-sphere.  We
also discuss three applications of this invariant to knot theory: studying the
Thurston norm of a link complement, the slice genus of a knot, and the
unknotting number of a knot. We emphasize the application to the
Thurston norm, and illustrate the theory in the case of the Conway link.
\end{abstract}

\maketitle
\section{Heegaard Floer homology of three-manifolds}

Floer homology was initially introduced by Floer to study questions in
Hamiltonian dynamics~\cite{FloerLag}. The basic set-up for his theory
involves a symplectic manifold $(M,\omega)$, and a pair of Lagrangian
submanifolds $L_0$ and $L_1$. His invariant, {\em Lagrangian Floer
  homology}, is the homology group of a chain complex generated freely
by intersection points between $L_0$ and $L_1$, endowed with a
differential which counts pseudo-holomorphic disks. This chain complex
arises from a suitable interpretation of the Morse complex in a certain 
infinite-dimensional setting.

Soon after formulating Lagrangian Floer homology, Floer
realized that his basic principles could also be used to construct a
three-manifold invariant, {\em instanton Floer homology}, closely
related to Donaldson's invariants for four-manifolds. In this version,
the basic set-up involves a closed, oriented three-manifold $Y$ (satisfying suitable
other topological restrictions on $Y$; for example, the theory is
defined when $Y$ has trivial integral first homology).  Again, one forms a
chain complex, but this time the generators are $SU(2)$
representations of the fundamental group of $Y$ (or some suitable
perturbation thereof), and the differentials count anti-self-dual
Yang-Mills connections on the product of $Y$ with the real line. This
invariant plays a crucial role in Donaldson's invariants for smooth
four-manifolds: for a four-manifold-with-boundary, the relative
Donaldson invariant is a homology class in the instanton Floer
homology groups of its boundary~\cite{DonaldsonFloer}.

In the present note, we will outline an adaptation of Lagrangian Floer
homology, {\em Heegaard Floer homology}, which gives rise to a closed
three-manifold invariant~\cite{HolDisk}, \cite{HolDiskTwo}. This
invariant also fits into a four-dimensional
framework~\cite{HolDiskFour}. There is a related invariant of smooth
four-manifolds, and indeed relative invariants for this four-manifold
invariant take values in the Heegaard Floer homology groups of its
boundary. 

A Heegaard diagram is a triple consisting of a closed, oriented
two-manifold $\Sigma$ of genus $g$, and a pair of $g$-tuples of
embedded, disjoint, homologically linearly independent curves
$\alphas=\{\alpha_1,...,\alpha_g\}$ and
$\betas=\{\beta_1,...,\beta_g\}$. A Heegaard diagram uniquely
specifies a three-manifold, obtained as a union of two genus $g$
handlebodies $U_\alpha$ and $U_\beta$. In $U_\alpha$, the curves
$\alpha_i$ bound disks, while in $U_\beta$, the curves $\beta_i$ bound
disks. We associate to this data a suitable version of Lagrangian
Floer homology.

Our ambient manifold in this case is the $g$-fold symmetric product of
$\Sigma$, the set of unordered $g$-tuples of points in $\Sigma$.  This
space inherits a natural complex structure from a complex structure
over $\Sigma$. Inside this manifold, there is a pair of
$g$-dimensional real tori, $\Ta=\alpha_1\times...\times\alpha_g$ and
$\Tb=\beta_1\times...\times\beta_g$. We fix also a reference point 
$$w\in\Sigma-\alpha_1-...-\alpha_g-\beta_1-...-\beta_g.$$
This gives
rise to a subvariety $V_w=\{w\}\times\Sym^{g-1}(\Sigma)\subset
\Sym^g(\Sigma)$. We consider the chain complex generated by
intersection points $\Ta\cap\Tb$. Concretely, an intersection point
of $\Ta$ and $\Tb$ corresponds to a permutation $\sigma$ in the symmetric
group on $g$ letters, together with a $g$-tuple of points $\x=(x_1,...,x_g)$
with $x_i\in\alpha_i\cap\beta_{\sigma(i)}$.

The differential again counts holomorphic disks; but some aspect of
the homotopy class of the disk is recorded. We make this precise
presently.  For fixed $\x,\y\in\Ta\cap\Tb$, let $\pi_2(\x,\y)$ denote the
space of homotopy classes of Whitney disks connecting $\x$ to $\y$,
i.e.  continuous maps of the unit disk $\CDisk\subset \C$
into $\Sym^g(\Sigma)$, mapping
the part of the boundary of $\CDisk$ with negative resp. positive real part to
$\Ta$ resp. $\Tb$, and mapping $i$ resp. $-i$ to $\x$ resp. $\y$.  The
algebraic intersection number of $\phi\in\pi_2(\x,\y)$ with the
subvariety $V_w$ determines a well-defined map $$n_w\colon
\pi_2(\x,\y)\longrightarrow \Z.$$ It is also useful to think of
the two-chain $\cald(\phi)$, which is gotten as a formal sum of
regions in $\Sigma-\alpha_1-...-\alpha_{g}-\beta_1-...-\beta_g$,
where a region is counted with multiplicity $n_{p}(\phi)$, 
where here $p\in\Sigma$ is any point in this region.
Given a Whitney disk, we can
consider its space of holomorphic representatives $\ModFlow(\phi)$,
using the induced complex structure on $\Sym^g(\Sigma)$.  If this
space is non-empty for all choices of almost-complex structure,
then the associated two-chain $\cald(\phi)$ has only non-negative
local multiplicities.
The group
$\R$ acts on $\ModFlow(\phi)$ by translation. The moduli space
$\ModFlow(\phi)$ has an expected dimension $\Mas(\phi)$, which is
obtained as the Fredholm index of the linearized $\DBar$-operator.
This quantity, the {\em Maslov index}, is denoted $\Mas(\phi)$.

It is sometimes necessary to perturb the holomorphic condition to
guarantee that moduli spaces are manifolds of the expected
dimension. It is useful (though slightly imprecise) to think of a
holomorphic disk in $\ModFlow(\phi)$ as a pair consisting of a
holomorphic surface $F$ with marked boundary, together with a degree
$g$ holomorphic projection map $\pi$ from $F$ to the standard disk,
and also a map $f$ from $F$ into $\Sigma$. Here, $f$ maps $\pi^{-1}$
of the subarc of the boundary of $\CDisk$ with negative resp. positive
real part into the subset $\alpha_1\cup...\cup\alpha_g$ resp
$\beta_1\cup...\cup\beta_g$.

We now consider the complex $\CFm(Y)$ which is the free
$\Z[U]$-module generated by $\Ta\cap\Tb$, with differential given by
\begin{equation}
\label{eq:DefD}
\partial \x = \sum_{\y\in\Ta\cap\Tb}
\sum_{\{\phi\in\pi_2(\x,\y)\big|\Mas(\phi)=1\}}
\#\left(\frac{\ModFlow(\phi)}{\R}\right)U^{n_{w}(\phi)}\y.
\end{equation}
In the case where $Y$ is an integral homology sphere, the above sum is
readily seen to be finite. (In the case where the first Betti number
is positive, some further constraints must be placed on the Heegaard
diagram.) With the help of Gromov's compactification of the space
pseudo-holomorphic curves~\cite{Gromov}, one can see that
$\partial^2=0$.

According to~\cite{HolDisk}, the homology groups $\HFm(Y)$ of
$\CFm(Y)$ are a topological invariant of $Y$. Indeed, the chain
homotopy type of $\CFm(Y)$ is a topological invariant, and, since
$\CFm(Y)$ is a module over $\Z[U]$, there are a number of other
associated constructions.  For example, we can form the  chain
complex $\CFinf(Y)$ obtained by inverting $U$, i.e.  a chain complex over
$\Z[U,U^{-1}]$, with differential as in Equation~\eqref{eq:DefD}.
The quotient of $\CFinf(Y)$ by $\CFm(Y)$ is a complex $\CFp(Y)$
which is often more convenient to work with. The corresponding
homology groups are denoted $\HFinf(Y)$ and $\HFp(Y)$ respectively.
Also, there is a chain complex $\CFa$ obtained by setting $U=0$;
explicitly, it is generated freely over $\Z$ by $\Ta\cap\Tb$, and
endowed with the differential
\[
\Da \x = \sum_{\y\in\Ta\cap\Tb}
\sum_{\{\phi\in\pi_2(\x,\y)\big|\Mas(\phi)=1, n_w(\phi)=0\}}
\#\left(\frac{\ModFlow(\phi)}{\R}\right)\y,
\]
and its homology (also a topological invariant of $Y$)
is denoted $\HFa(Y)$.

The invariants $\HFm(Y)$, $\HFinf(Y)$, and $\HFp(Y)$, together with
the exact sequence connecting them, are crucial ingredients in the
construction of a Heegaard Floer invariant $\Phi$ for closed, smooth
four-manifolds.  We will say only little more about this invariant
here, referring the reader to~\cite{HolDiskFour}
for its construction.

\section{Heegaard Floer homology of knots}
Heegaard Floer homology for three-manifolds has a refinement to an
invariant for null-homologous knots in a three-manifold, as defined
in~\cite{Knots}, and also independently by Rasmussen
in~\cite{RasmussenThesis}.

A knot $K$ in a three-manifold $Y$ is specified by a Heegaard diagram
$(\Sigma,\alphas,\betas)$ for $Y$, together with a pair $w$ and $z$ of
basepoints in $\Sigma$. The knot $K$ is given as follows. Connect $w$
and $z$ by an arc $\xi$ in $\Sigma-\alpha_1-...-\alpha_g$ and an arc
$\eta$ in $\Sigma-\beta_1-...-\beta_g$. The arcs $\xi$ and $\eta$ are
then pushed into $U_\alpha$ and $U_\beta$ respectively, so that they
both meet $\Sigma$ only at $w$ and $z$, giving new
arcs $\xi'$ and $\eta'$. Our knot $K$, then, is given
by $\xi'-\eta'$.  For simplicity, we consider here the case where the
ambient manifold $Y$ is the three-sphere $S^3$.

The new basepoint $z$ gives the Heegaard Floer complex a filtration.
Specifically, we can construct a map
$$F\colon \Ta\cap\Tb\longrightarrow \Z$$
by 
\begin{equation}
\label{eq:DefFilt}
F(\x)-F(\y)=n_z(\phi)-n_w(\phi),
\end{equation}
where $\phi\in\pi_2(\x,\y)$. It is easy to see that this quantity
is independent of the choice of $\phi$, depending only on $\x$ and $\y$.
Moreover, if $\y$ appears with non-zero multiplicity
in $\Da(\x)$, then $F(\x)\geq F(\y)$. This follows from the fact that
there is a pseudo-holomorphic disks $\phi\in\pi_2(\x,\y)$ with $n_w(\phi)=0$,
and also $n_z(\phi)\geq 0$, since a pseudo-holomorphic disks meets
the subvariety $V_z$ with non-negative intersection number.

Equation~\eqref{eq:DefFilt} defines $F$ uniquely up to an overall translation.
This indeterminacy will be removed presently.

The filtered chain homotopy type of this filtered chain complex is an invariant
of the knot $K$.  For example, the homology of the associated graded object,
the {\em knot Floer homology} is an invariant of $K\subset S^3$, defined by
$$\HFKa(S^3,K)=\bigoplus_{s\in\Z}\HFKa(S^3,K,s),$$
where $\HFKa(S^3,K,s)$ is the homology group of the chain complex
generated by intersection points $\x\in\Ta\cap\Tb$ with $F(\x)=s$, endowed
with differential 
\[
\partial \x = \sum_{\y\in\Ta\cap\Tb}
\sum_{\left\{\phi\in\pi_2(\x,\y)\big|
\begin{tiny}\begin{array}{c}
\Mas(\phi)=1, \\
n_w(\phi)=n_z(\phi)=0 
\end{array}\end{tiny}\right\}}
\#\left(\frac{\ModFlow(\phi)}{\R}\right)\y.
\]

The graded Euler characteristic of this theory is the Alexander polynomial of $K$,
in the sense that
\begin{equation}
  \label{eq:KnotEuler}
\Delta_K(T)=\sum_{s\in\Z} \chi(\HFKa_*(K,s)) \cm T^{s}.
\end{equation}
This
formula can be used to pin down the additive indeterminacy of $F$: we
require that $F$ be chosen so that the graded Euler characteristic is
the symmetrized Alexander polynomial. In fact, this symmetry has a stronger
formulation, as a relatively graded isomorphism
\begin{equation}
  \label{eq:KnotSymmetry}
  \HFKa_*(K,s)\cong \HFKa_*(K,-s).
\end{equation}

\section{Heegaard Floer homology for links}

Heegaard Floer homology groups of knots can be generalized to the case
of links in $S^3$. For an $\ell$-component link, we
consider a Heegaard diagram with genus $g$ Heegaard surface, and two
$(g+\ell-1)$-tuples attaching circles
$\alphas=\{\alpha_1,...,\alpha_{g+\ell-1}\}$ and
$\betas=\{\beta_1,...,\beta_{g+\ell-1}\}$. We require
$\{\alpha_1,...,\alpha_{g+\ell-1}\}$ to be disjoint and embedded, and
to span a $g$-dimensional lattice in $H_1(\Sigma;\Z)$. The same is
required of the $\{\beta_1,...,\beta_{g+\ell-1}\}$. Clearly,
$\Sigma-\alpha_1-...-\alpha_{g+\ell-1}$ consists of $\ell$ components
$A_1,...,A_\ell$.  Similarly, $\Sigma-\beta_1-...-\beta_{g+\ell-1}$
consists of $\ell$ components $B_1,...,B_\ell$. We assume that this
Heegaard diagram has the special property that $A_i\cap B_i$ is
non-empty. Indeed, for each $i=1,...,\ell$, we choose basepoints
$w_i$ and $z_i$ to lie inside $A_i\cap B_i$. We call the collection of data
$(\Sigma,\alphas,\betas,\{w_1,...,w_\ell\},\{z_1,...,z_\ell\})$ a
{\em $2\ell$-pointed Heegaard diagram}.

A link can now be constructed in the following manner. Connect $w_i$
and $z_i$ by an arc $\xi_i$ in $A_i$ and an arc $\eta_i$ in $B_i$.
Again, the arc $\xi_i$ resp. $\eta_i$ is pushed into $U_\alpha$ resp.
$U_\beta$ to give rise to a pair of arcs $\xi_i'$ and $\eta_i'$.  The
link $L$ is given by $\cup_{i=1}^\ell \xi_i'-\eta_i'$.  For a
$2\ell$-pointed Heegaard diagram for $S^3$
$(\Sigma,\alphas,\betas,\{w_1,...,w_\ell\}, \{z_1,...z_\ell\})$, if
$L$ is the link obtained in this manner, we say that
the Heegaard diagram is {\em compatible} with the link $L$.

We will need to make an additional assumption on the Heegaard diagram.
A {\em periodic domain} is a two-chain in $\Sigma$ of the form 
$$\sum c_i (A_i-B_i),$$ 
where $c_i\in\Z$. Our assumption is that all non-zero periodic domains
have some positive and some negative local multiplicities $c_i$.
This assumption on the pointed Heegaard diagram is called {\em admissibility}.

Let $L\subset S^3$ be an $\ell$-component link,
suppose that $L$ is embedded so that the restriction of the height function
to $L$ has $b$ local maxima, then 
we can construct a compatible $2\ell$-pointed heegaard diagram with Heegaard
genus $g=b-\ell$.

\begin{figure}[ht]
\mbox{\vbox{\epsfbox{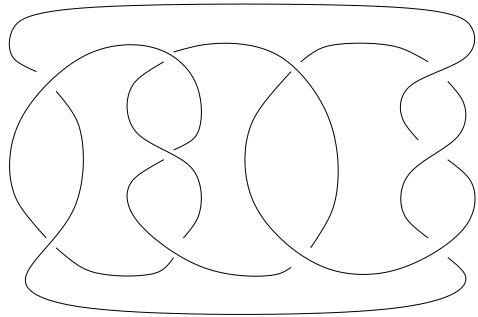}}}
\caption{\label{fig:Conway}
{\bf{The Conway link.}}}
\end{figure}

For example, consider the two-component ``Conway link'' pictured in
Figure~\ref{fig:Conway}. This is the $(2,-3,-2,3)$ pretzel link,
also known as L10n59 in Thistlethwaite's link table~\cite{Thistlethwaite}.

\begin{figure}[ht]
\mbox{\vbox{\epsfbox{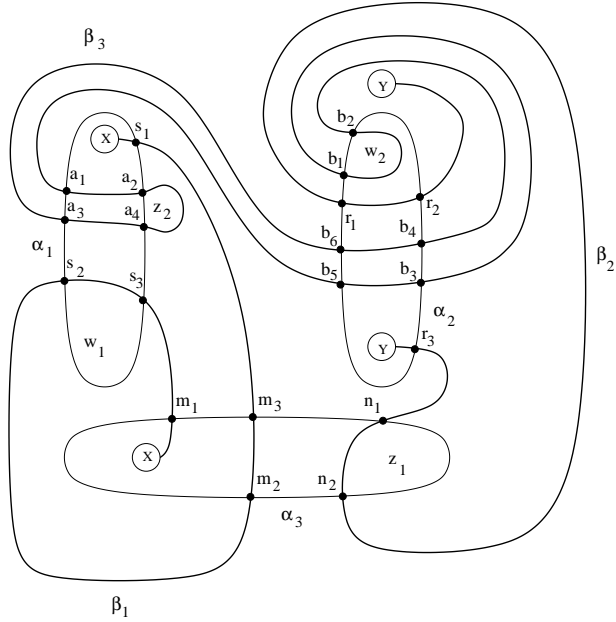}}}
\caption{\label{fig:ConwayDiag}
{\bf{Pointed Heegaard diagram for the Conway link.}}
This picture takes place on the genus two surface obtained
by identifying the two disks labeled by $X$ and the two disks
labeled by $Y$.}
\end{figure}

For this link, $b=4$, and hence we can draw it on a surface of genus $g=2$,
as illustrated in Figure~\ref{fig:ConwayDiag}.
It is straightforward to verify that the space of periodic domains
is one-dimensional; drawing a picture of this generator, it is
also straightforward to see that the diagram is admissible.

Now, we work inside $\Sym^{g+\ell-1}(\Sigma)$, relative to the tori
$\Ta=\alpha_1\times...\times\alpha_{g+\ell-1}$ and
$\Tb=\beta_1\times...\times\beta_{g+\ell-1}$, and consider
intersection points of $\Ta$ and $\Tb$; i.e. $g+\ell-1$-tuples of
points $(x_1,...,x_{g+\ell-1})$ with $x_i\in
\alpha_i\cap\beta_{\sigma(i)}$ for some permutation $\sigma$ in the
symmetric group on $g+\ell-1$ letters. We then form the chain complex
$\CFLa(S^3,L)$ generated freely by these intersection points.

For example, for the figure illustrated in Figure~\ref{fig:ConwayDiag},
the curves $\alpha_i$ and $\beta_j$ intersect according to the pattern
\vskip.2cm
\begin{center}
\begin{tabular}{c|ccc}
$\cap$ &  $\alpha_1$ & $\alpha_2$ & $\alpha_3$  \\
        \hline
$\beta_1$ & $\{s_1,s_2,s_3\}$ & $\emptyset$ & $\{m_1,m_2,m_3\}$\\
$\beta_2$  & $\emptyset$ &  $\{r_1,r_2,r_3\}$ & $\{n_1,n_2\}$ \\
$\beta_3$ & $\{a_1,...,a_4\}$ & $\{b_1,...,b_6\}$ & $\emptyset$
\end{tabular}
\end{center}
\vskip.2cm
Now, there are exactly two permutations of $\{1,2,3\}$ for which
$\alpha_i\cap\beta_{\sigma(i)}$ is non-trivial for all $i$. This gives
two types of intersection points of $\Ta\cap\Tb$, namely,
$a_i\times m_j \times r_k$ (with $i=1,...,4$,
$j=1,...,3$, $k=1,...,3$) and also $b_i\times n_j\times s_k$
(with $i=1,...,6$, $j=1,2$, $k=1,2,3$). This gives
a chain complex with a total of $72$ generators.

The complex $\CFLa$ has a grading, the {\em Maslov grading},
which is specified up to overall translation by the convention
$$\gr(\x)-\gr(\y)=\Mas(\phi)-2\sum_{i=1}^{\ell} n_{w_i}(\phi),$$
where
$\phi\in\pi_2(\x,\y)$ is any Whitney disk connecting $\x$ and $\y$.
The parity of the Maslov grading depends
on the local sign of the intersection number of $\Ta$ and $\Tb$
at $\x$. 

But $\CFLa$ has an additional grading, the $\HH$-grading. To define
this, we associate to each $\phi\in\pi_2(\x,\y)$ the pair of
vectors
\begin{eqnarray*}
n_{\w}(\phi)=(n_{w_1}(\phi),...,n_{w_\ell}(\phi))
&{\text{and}}&
n_{\z}(\phi)=(n_{z_1}(\phi),...,n_{z_\ell}(\phi)).
\end{eqnarray*}
We have a function $F\colon \Ta\cap\Tb\longrightarrow
\Z^\ell\cong H_1(S^3-L;\Z)$ 
(where the latter identification
is given by the meridians of the link $L$)
specified uniquely up to translation by the formula
$$F(\x)-F(\y)=n_{\z}(\phi)-n_{\w}(\phi),$$
where $\phi$ is any choice of homotopy class in $\pi_2(\x,\y)$.

Endow $\CFLa(S^3,L)$ with the differential
\[
\partial \x = \sum_{\y\in\Ta\cap\Tb}
\sum_{\left\{\phi\in\pi_2(\x,\y)\big|
\begin{tiny}\begin{array}{c}
\Mas(\phi)=1, \\
n_{\w}(\phi)=n_{\z}(\phi)=0 
\end{array}\end{tiny}\right\}}
\#\left(\frac{\ModFlow(\phi)}{\R}\right)\y.
\]
It is easy to see that this differential drops Maslov grading by one.
Moreover, the complex naturally splits into summands indexed by elements 
of $\Z^\ell\cong H_1(S^3-L;\Z)$ specified by the function $F$.
We find it natural to think of these summands, in fact,
as indexed by the affine lattice
$\HH\subset H_1(S^3-L;\R)$  over
$H_1(S^3-L;\Z)$, given by elements $$\sum_{i=1}^\ell a_i\cm [\mu_i],$$
where $a_i\in\Q$ satisfies the property that $$2a_i+\lk(L_i,L-L_i)$$
is an even integer. The translational ambiguity of the 
map is then pinned down by the following generalization of
Equation~\eqref{eq:KnotSymmetry}:
\begin{equation}
  \label{eq:Symmetry}
  \HFLa_*(\orL,h)\cong \HFLa_*(\orL,-h).
\end{equation}

In practice, it is easy to calculate the difference in $F$ for any two
intersection of $\Ta$ and $\Tb$ which have the same type (i.e. same
pattern of intersection $\alpha_i\cap\beta_{\sigma(i)}$). To this end,
it suffices to find for each pair of intersection points $x, x'\in
\alpha_i\cap\beta_j$, a disk (or more generally a compact surface with
a single boundary component) in $\Sigma$ which meets $\alpha_i$ along
one arc in its boundary and $\beta_j$ along the complementary arc,
carrying the intersection points of the closures of the arcs to $x$
and $x'$. We then define the ``relative difference'' of
$x$ and $x'$, $F^{i,j}(x)-F^{i,j}(x')$, to be
$n_{\z}-n_{\w}$ for this disk (or surface). It is easy to see
then that if $\x$ and $\y\in\Ta\cap\Tb$ are two intersection points
with the same type (as specified by $\sigma$), then
$$F(\x)-F(\y)=\sum_{i=1}^{g+\ell-1} F^{i,\sigma(i)}(x_i)-
F^{i,\sigma(i)}(y_i),$$
where $\x=(x_1,...,x_{g+\ell-1})$ and
$y=(y_1,...,y_{g+\ell-1})$.  This determines $F(\x)-F(\y)$ for $\x$
and $\y$ of the same type.  Different types can then be compared by
choosing homotopy classes connecting them (and in suitable
circumstances, the translational ambiguity can be removed using
Equation~\eqref{eq:Symmetry}).

\begin{figure}[ht]
\mbox{\vbox{\epsfbox{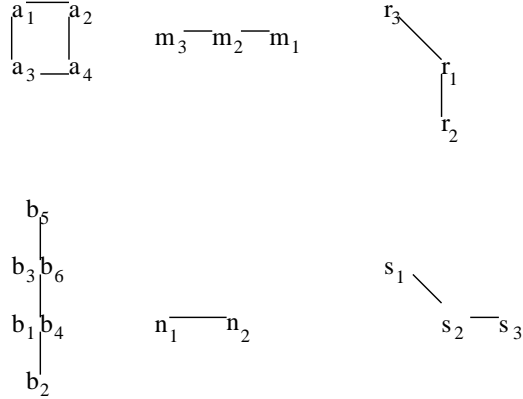}}}
\caption{\label{fig:ConwayGens}
{\bf{Generators for $\HFLa$ of the Conway link.}}
We illustrate the relative differences of the various intersection
points of $\alpha_i$ and $\beta_j$. A horizontal resp. vertical segment denotes
two intersection points whose relative difference is one in the
first resp. second component; e.g. there is a disk in Figure~\ref{fig:ConwayDiag}
$\phi$ from $a_2$ to $a_1$ with $n_{\z}-n_{\w}$ given by $(1,0)$,
while there is one from $b_5$ to $b_3$ with relative difference
given by $(0,1)$. Finally, for the diagonal edges, we have a disk
from $r_3$ to $r_1$ with relative difference $(-1,1)$.}
\end{figure}

We  display the relative differences for the various
intersection points for the diagram from Figure~\ref{fig:ConwayGens}.

\begin{figure}[ht]
\mbox{\vbox{\epsfbox{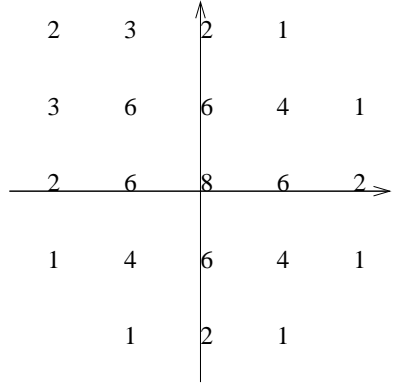}}}
\caption{\label{fig:RkCFL}
{\bf{Generators for $\CFLa$ of the Conway link.}}
The $72$ generators of $\CFLa$ coming from the diagram
in Figure~\ref{fig:ConwayDiag} are separated into different .}
\end{figure}

It is now straightforward to verify that the ranks of the chain 
groups in each value of $F$ is given as in Figure~\ref{fig:RkCFL}.
It is  more challenging to calculate the homology
groups of $\CFLa$.

Some aspects are immediate. For example, it follows by glancing at
Figure~\ref{fig:RkCFL}, and comparing Equation~\eqref{eq:Symmetry}
that the homology in $\HH$-grading $(-2,2)$ is trivial (as
there are no generators in the $\HH$-grading $(2,-2)$), and that in
$\HH$-gradings $(-1,2)$ and $(-2,1)$ the groups $\HFLa$ have rank one.
This already suffices to determine the convex hull of $h\in \HH$
for which $\HFLa(L,h)$ is non-trivial, as required for the application to 
the Thurston norm below (see esp. Equation~\eqref{eq:DetectsThurstonNorm}).

Also, the calculation of $\HFLa(L,(x,y))$ with $(x,y)\in \{(0,\pm 2),
(\pm 2,0)\}$ follows from the fact that for each of these $\HH$-gradings,
every generator has the same Maslov grading.

\begin{figure}[ht]
\mbox{\vbox{\epsfbox{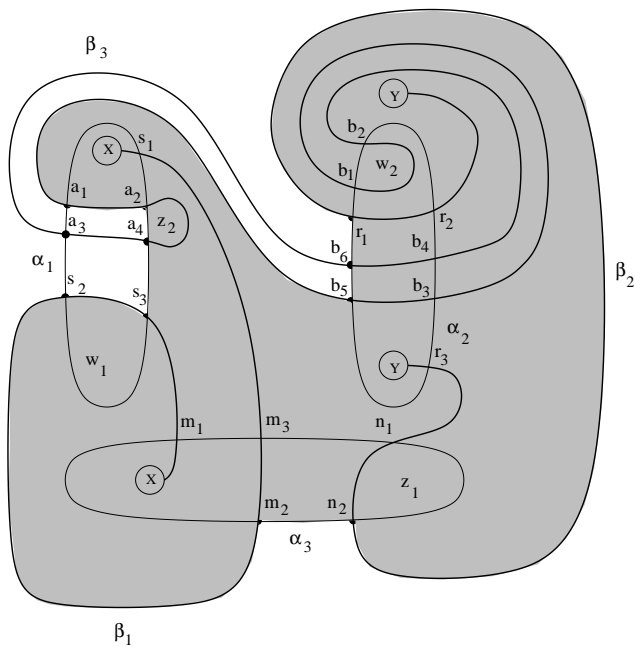}}}
\caption{\label{fig:AHexagon}
{\bf{A flowline.}}
The complement of the shaded region gives a hexagon
connecting $b_5\times n_1\times s_3$ to 
$a_2\times m_2\times r_1$.}
\end{figure}

Next, consider the part in $\HH$-grading $(1,1)$. There are four generators
$$\{
\begin{array}{llll}
a_1\times m_1\times r_1, &a_2\times m_2\times r_1, &
a_4\times r_3\times m_1,& b_5\times n_1\times s_3
\end{array}\}.$$
For the case where $\x=b_5\times n_1\times s_3$ and
$\y\in\{a_1\times m_1\times r_1, 
a_2\times m_2\times r_1\}$, there is a homotopy class
$\phi\in\pi_2(\x,\y)$ whose associated two-chain $\cald(\phi)$
is a hexagon.
For the case where $\y=a_2\times m_2\times r_1$,
we illustrate this in Figure~\ref{fig:AHexagon}.
A hexagon gives rise to a flow-line connecting $\x$ to $\y$. 
To this end, we think of a holomorphic disk in $\Sym^3(\Sigma)$ as
a branched triple-cover $F$ of the standard disk, together with a map 
of $F$ into $\Sigma$. The given hexagonal domain in $\Sigma$ can
be realized as a branched triple-cover of the disk $\CDisk$. 
Moreover, any other domain connecting $\x$ to $\y$ has negative local
multiplicity somewhere. Hence, we 
have that
$$\partial b_5\times n_1\times s_3
= a_1\times m_1\times r_1 + 
a_2\times m_2\times r_1.$$
It can also be seen that $$\gr(b_5\times n_1\times s_3)=
\gr(a_4\times m_1\times r_3),$$
but 
there are no non-negative domains from $a_4\times m_1\times r_3$
to either of $\{a_1\times m_1\times r_1, 
a_2\times m_2\times r_1\}$. It follows at once that $\HFLa(L,(1,1))$
has rank two.

With some additional work, one can verify that
all the link Floer homology groups of the Conway link 
are as displayed below in Figure~\ref{fig:RkHFL}.
\begin{figure}[ht]
\mbox{\vbox{\epsfbox{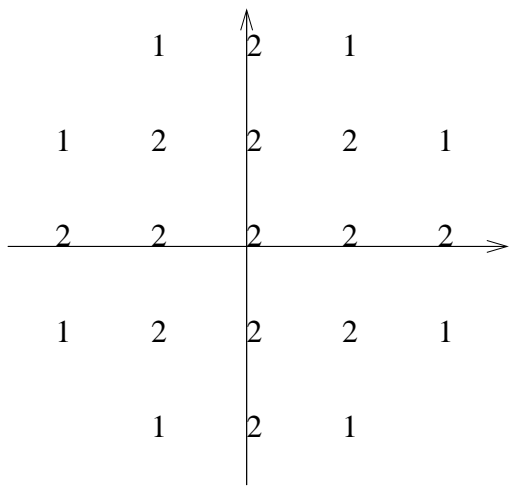}}}
\caption{\label{fig:RkHFL}
{\bf{Ranks of $\HFLa$ for the Conway link.}}
The ranks are displayed, along with their $\HH$ grading, thought
of as coordinates in the plane.}
\end{figure}

\section{Basic properties}

Perhaps the single most fundamental property of Heegaard Floer
homology is that it satisfies an exact triangle for surgeries.
More precisely, a {\em triad} of three-manifolds $Y_1$, $Y_2$, $Y_3$
is a cyclically ordered triple of three-manifolds obtained as follows.
Let $M$ be a three-manifold with torus boundary, and fix three 
simple, closed curves in its boundary $\gamma_1$, $\gamma_2$, $\gamma_3$
any two of which intersect transversally in one point, and ordered
so that there are orientations on the three curves so that
$$
\#(\gamma_1\cap \gamma_2)=\#(\gamma_2 \cap \gamma_3)
=\#(\gamma_3\cap \gamma_1) = -1.
$$
We let $Y_i$ be the three-manifold obtained by Dehn filling $M$ along
the curve $\gamma_i$.

\begin{theorem}
  \label{thm:ExactTriangle}
  Let $Y_1$, $Y_2$, and $Y_3$ be a triad of three-manifolds. Then, there 
  is an long exact sequence of the form
  $$ 
  \begin{CD}
    ...@>>>\HFa(Y_1)@>>>\HFa(Y_2) @>>> \HFa(Y_3)@>>> ...
  \end{CD}.
  $$
\end{theorem}

The maps in the exact triangle are induced by the three natural
two-handle cobordisms connecting $Y_i$ and $Y_{i+1}$ (where we $i$ as
an integer modulo $3$), in a manner made precise in~\cite{HolDiskFour}.

The above surgery exact sequence is similar
to an exact sequence
established by Floer for instanton Floer homology (only using a
restricted class of triads)~\cite{FloerTriangles},
\cite{BraamDonaldson}. An analogous exact sequence has been
established for Floer homology of Seiberg-Witten monopoles,
see~\cite{KMOSz}. There are also related exact sequences in
symplectic geometry, cf.~\cite{Seidel}.

Note that there are other variants of the exact triangle, and indeed,
there are certain other related calculational techniques for Heegaard
Floer homology. For example, for a knot $K$ in an integral homology
sphere $Y$, the filtered chain homotopy type of the induced knot
invariant can be used to calculate the Floer homology groups of
arbitrary surgeries on $K$,~\cite{RatSurg}.

\section{Three applications}

Heegaard Floer homology is particularly well suited to problems in
knot-theory and three-manifold topology which can be formulated in
terms of the existence of four-dimensional cobordisms. We focus here
on a few concrete problems which can be formulated for knots and links
in the three-sphere. For some other applications, see~\cite{NoteLens},
\cite{HolDiskContact}, \cite{LiscaStipsicz}, \cite{RasmussenLens}.

\subsection{Thurston norm}
Let $K\subset S^3$ be a knot. The {\em Seifert genus} of $K$, denoted
$g(K)$, is the minimal genus of any embedded surface $F\subset S^3$
with boundary $K$. Clearly, if $g(K)=0$, then $K$ is the unknot.

According to~\cite{GenusBounds}, knot Floer homology detects  the
Seifert genus of a knot, by the property that
\begin{equation}
  \label{eq:DetectsGenus}
  g(K)=\max\{s\big| \HFKa(K,s)\neq 0\}.
\end{equation}

There is a natural generalization of the knot genus and indeed
of Equation~\eqref{eq:DetectsGenus}. This is best 
formulated in terms of Thurston's norm
on second homology.

Recall that if $F$ is a compact, oriented, but possibly disconnected
surface-with-boundary $F=\bigcup_{i=1}^n F_i$, its {\em complexity}
is given by
$$\chi_-(F)=\sum_{\{F_i\big| \chi(F_i)\leq 0\}} -\chi(F_i).$$
Given
any homology class $h\in H_2(S^3,L)$, it is easy to see that there is
a compact, oriented surface-with-boundary embedded in $S^3-\nbd{K}$
representing $h$.  Consider the function from $H^1(S^3-L;\Z)$ to the
integers defined by
$$x(h)=\min_{\{F\hookrightarrow S^3-\nbd{K}\big|[F]=\PD[h]\}} \chi_-(F).$$
Indeed, according to Thurston~\cite{Thurston},
this function $x$ satisfies an inequality
$x(h_1+h_2)\leq x(h_1)+x(h_2)$,
and it is linear on rays, i.e. given $h\in H_2(S^3,L)$
and a non-negative integer $n$, we have that $x(n\cm h)=n x(h)$.
Thus, $x$ can  
be naturally extended to a semi-norm
on $H^1(S^3-L;\R)$, the {\em Thurston
semi-norm}. In fact, this semi-norm is uniquely specified by 
its unit ball
$$B_x = \{h\in H^1(S^3-L;\R)\big| x(h)\leq 1\},$$
which is a polytope $H^1(S^3-L;\R)$ whose vertices lie at lattice
points in $H^1(S^3-L;\Z)$.

Equation~\eqref{eq:DetectsGenus} can now be generalized as follows.
A {\em trivial component} of a link $L$ is a component $K\subset L$ which is unknotted
an geometrically unlinked from the complement $L-K$. Suppose that $L$ is a link
with no trivial components. Then, 
given $s\in H^1(S^3-L)$, 
\begin{equation}
  \label{eq:DetectsThurstonNorm}
  x(h)+\sum_{i=1}^\ell |\langle \mu_i, h\rangle| =
2\cm \max_{\{s\in \HH\big| \HFLa(L,s)\neq 0\}}\langle h, s\rangle,
\end{equation}
where here $\langle,\rangle$ denotes the Kronecker pairing of
$H_1(S^3-L;\R)$ with $H^1(S^3-L;\R)$. 

This formula can be thought of more
geometrically from the following point of view. Consider the dual norm
$x^*\colon H_1(S^3-L;\R)\longrightarrow \R$ given by
$$x^*(s)=\max_{\{h\in B_x\}} \langle s,h\rangle.$$
The unit ball $B_{x*}$ is a (possibly degenerate) polytope in $H_1(S^3-L;\R)$
called the {\em dual Thurston polytope}.
Equation~\eqref{eq:DetectsThurstonNorm} states that for a link  $L$  with
no trivial components,
if we take the convex
hull of the set of $s\in \HH$ with $\HFLa(L,s)$, and rescale that polytope by a factor
of two, then we obtain the sum of 
the dual Thurston polytope with the symmetric hypercube in $H_1(S^3-L;\R)$ 
with edge-length two.

Of course, the Thurston norm can be defined for closed
three-manifolds, as well.  In fact, a result analogous to
Equation~\eqref{eq:DetectsThurstonNorm} can be proved for closed
three-manifolds $Y$, instead of link complements. An analogous result
has been shown to hold for Seiberg-Witten monopole Floer
homology~\cite{KMscalar} (but at present there is no analogue of knot
and link Floer homology in gauge-theoretic terms).
 
Although the statement of Equation~\eqref{eq:DetectsGenus} does not
explicitly involve any four-dimensional theory, the proof of this
result does use the full force of Heegaard Floer homology, combined
with Gabai's theory of sutured manifolds, and recent results in
symplectic geometry.  Specifically, according to a combination of
theorems of Gabai~\cite{Gabai}, \cite{GabaiKnots},
Eliashberg-Thurston~\cite{EliashbergThurston}, and a result of
Eliashberg~\cite{Eliashberg} and independently Etnyre~\cite{Etnyre},
if $K\subset S^3$ is a knot of genus $g$, then its zero-surgery
$S^3_0(K)$ separates a symplectic manifold. Non-vanishing theorems for
the Heegaard Floer invariant $\Phi$ for symplectic four-manifolds,
which in turn are built on the symplectic Lefshetz pencils of
Donaldson~\cite{Donaldson}, then give a non-vanishing result for the
Heegaard Floer $S^3_0(K)$ from which Equation~\eqref{eq:DetectsGenus}
follows.

These results can be further generalized to give
Equation~\eqref{eq:DetectsThurstonNorm}.  Specifically, an
$n$-component link in $S^3$ naturally gives rise to a connected knot
in the $(n-1)$-fold connected sum of $S^2\times S^1$. A genus bound
analogous to Equation~\eqref{eq:DetectsGenus} has been shown by Ni
in~\cite{Ni}, which in effect establishes
Equation~\eqref{eq:DetectsThurstonNorm}, in the case where $h$ is one
of the $2^\ell$ cohomology classes whose evaluation on each meridian
for $L$ has absolute value equal to one.
Equation~\eqref{eq:DetectsThurstonNorm} then follows from the manner
in which the Thurston norm and link Floer homology transform under
cabling, see also~\cite{Cables}.

\begin{figure}[ht]
\mbox{\vbox{\epsfbox{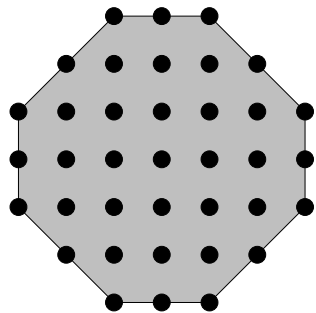}}}
\caption{\label{fig:ConwayThurston}
{\bf{Dual Thurston polytope for the Conway link.}}
Lattice points in $H_1(S^3-L;\Z)$ are indicated by solid
circles.}
\end{figure}

As an illustration of Equation~\eqref{eq:DetectsGenus}, consider
the Conway link from Figure~\ref{fig:Conway}. According to the calculations
displayed in Figure~\ref{fig:RkHFL}, together with this equation,
we conclude that the dual Thurston polytope of the Conway link is as illustrated
in Figure~\ref{fig:ConwayThurston}.

Figure~\ref{fig:ConwayThurston} suggests that there are surfaces $F_1$
the complement of the Conway link $L$, with $\partial F_1$ consisting
of a longitude belonging to a first component of $L$, and some number
of copies of the meridian of the second component, and also with
$\chi_-(F_1)=5$.  In fact, such a surface can be easily obtained by
puncturing a genus one Seifert surface for one of the trefoil
components in two additional points. A similar surface can be found
for the other component of $L$. 

Note that a verification of the dual Thurston polytope for the Conway
link can be easily obtained by more classical methods
(cf.~\cite{McMullen}); however, the computation given here is fairly
easy (and hopefully illustrates the theory).

\subsection{Slice genus}

A {\em slice surface} for a knot $K$ is a smoothly embedded
surface-with-boundary $F\subset B^4$ which meets $S^3$ along its
boundary, which is $K$.  The {\em slice genus} of a knot $g_*(K)$ is
the minimal genus of any slice surface for $K$.  Heegaard Floer
homology can be used to give information about this quantity, as
follows.

Recall that knot Floer homology is the homology of an associated
graded object which is induced by the filtration of a chain complex
which calculates $\HFa(S^3)\cong \Z$. But the entire filtered chain
homotopy type of the complex is a knot invariant. Denote the sequence
of subcomplexes $F_i\subset F_{i+1}$, so that for all sufficiently
small integers $i$, $F_i=0$, while for all sufficiently large
integers, $F_i=\CFa(S^3)$.  There is another natural invariant which
can be associated to a knot, which is the minimal $i$ for which the
map $H_*(F_i) \longrightarrow \HFa(S^3)$ is non-trivial. According
to~\cite{4BallGenus} and independently~\cite{RasmussenThesis},
$$|\tau(K)|\leq g^*(K).$$
Intriguingly,
Rasmussen~\cite{RasmussenSlice} has shown that a very similar
algebraic construction on Khovanov's homology~\cite{Khovanov},
\cite{EunSooLee}, can be used to define a similar (but entirely combinatorial) 
numerical invariant
$s(K)$.  Although both $\tau(K)$ and $s(K)$ share many formal
properties (and hence agree on many knots), Hedden and Ording have recently
shown~\cite{HeddenSTau} that these two invariants are in fact distinct.
Their examples are certain twisted Whitehead doubles of the trefoil.

\subsection{Unknotting numbers}
The {\em unknotting number} $u(K)$ is the minimal number of crossing
changes required to transform $K$ into an unknot.  An $n$-step
unknotting of a knot $K$ in effect gives an immersed disk in $B^4$
with $n$ double-points. Resolving these double-points, we obtain a
slice surface for $K$ with genus $n$.  This observation immediately
verifies the inequality $$g^*(K)\leq u(K).$$

For some classes of knots, these two quantities are equal. For
example, for the $(p,q)$ torus knot, $g^*(K)=u(K)=(p-1)(q-1)/2$.  This
was first shown by Kronheimer and Mrowka in~\cite{KMMilnor} (though it
has alternative proofs now using either $\tau$~\cite{4BallGenus} or the
Khovanov-Rasmussen invariant $s$~\cite{RasmussenSlice}). 

But there are Floer-theoretic bounds on $u(K)$ which go beyond the slice genus,
cf.~\cite{UnknotOne}, \cite{Owens}.

Suppose that $K$ has $u(K)=1$. Then, Montesinos
observed~\cite{Montesinos} that the branched double-cover of $S^3$
with branching locus $K$, denoted $\Sigma(K)$, can be realized as $\pm
d/2$-surgery on a different knot $C\subset S^3$, where here
$d=|\Delta_K(-1)|$. Obstructions to this can sometimes be given using
Heegaard Floer homology.

To do this in a useful manner, we must understand first
$\HFa(\Sigma(K))$.  For some knots, this is a straightforward matter.
For example, when $K$ is a knot which admits an alternating
projection, an easy induction using Theorem~\ref{thm:ExactTriangle}
shows that $\HFa(\Sigma(K))$ is a free $\Z$-module of rank
$|\Delta_K(-1)|$. This means that the Heegaard Floer homology groups
of these three-manifolds is as simple as possible. For any rational
homology three-sphere $Y$ (i.e. closed three-manifold with
$H_1(Y;\Q)=0$), the Euler characteristic of $\HFa(Y)$ is
$|H_1(Y;\Z)|$, the number of elements in $H_1(Y;\Z)$.  A rational
homology three-sphere whose homology group $\HFa(Y)$ is a free module
of rank $|H_1(Y;\Z)|$ is called an {\em $L$-space}.  Thus, if $K$ is a
knot with alternating projection, then $\Sigma(K)$ is an $L$-space.

\begin{figure}[ht]
\mbox{\vbox{\epsfbox{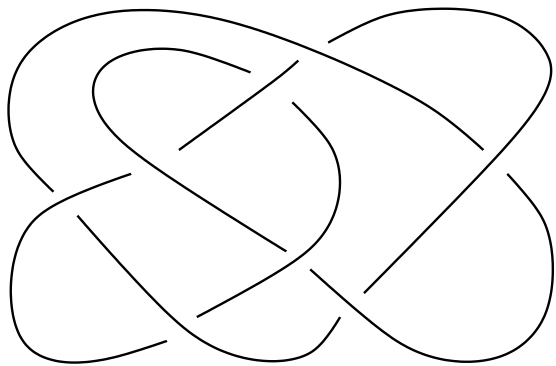}}}
\caption{\label{fig:8s10}
{\bf{A knot with $u=2$.}}}
\end{figure}

There are obstructions to realizing an $L$-space as surgery on a knot
in $S^3$, These obstructions are phrased in terms of an additional
$\Q$-grading on the Heegaard Floer homology~\cite{UnknotOne},
analogous to an invariant defined by Fr{\o}yshov~\cite{Froyshov} in
the context of Seiberg-Witten theory.  Moreover, this $\Q$-grading can
be explicitly calculated for $\Sigma(K)$ for an alternating knot $K$
from its Goeritz matrix.  Rather than stating these results precisely,
we content ourself here with including a picture of an eight-crossing
alternating knot ($8_{10}$, see Figure~\ref{fig:8s10}) whose unknotting number can be shown to
equal two via these (and presently, no other known) techniques.

Combining these obstructions with recent work of Gordon and
Luecke~\cite{GLUnknotOne}, one can classify all knots with $10$
and fewer crossings which have unknotting number equal to one. Indeed,
a different application of Heegaard Floer homology along similar lines
discovered by Owens~\cite{Owens} can be used to complete the
unknotting number table for prime knots with nine or fewer crossings.

\bibliographystyle{plain}
\bibliography{biblio}

\end{document}